\documentclass[reqno]{amsart}
\usepackage{url}
\usepackage{amssymb,amsthm,amsfonts,amstext,fancyhdr,amsmath,mathrsfs,epsfig,color}
\usepackage{amsmath}
\usepackage{mathptmx}       
\usepackage{helvet}         
\usepackage{courier}        
\usepackage{graphicx}
\usepackage{enumerate}
\usepackage{a4wide}
\usepackage[active]{srcltx}
\usepackage[english]{babel}
\usepackage[latin1]{inputenc}
\usepackage{color}

\newtheorem{theorem}{Theorem}

\newtheorem{corollary}[theorem]{Corollary}

\newtheorem{definition}[theorem]{Definition}

\newtheorem{proposition}[theorem]{Proposition}
\newtheorem{remark}[theorem]{Remark}

\begin{document}

\title{A formula for the Entropy of the Convolution of Gibbs probabilities on the circle}
\author{ Artur O. Lopes \, } 
\date{\today}

\maketitle

\centerline{Inst. de Matem\'atica, UFRGS - Porto Alegre, Brasil}
\bigskip


\vspace*{0.5 cm}
%
%
\bigskip

\begin{abstract} Consider the transformation $T:S^1 \to S^1$, such that $T(x)=2\, x$ (mod 1), and where $S^1$ is the unitary circle.
Suppose $J:S^1 \to \mathbb{R}$ is H\"older continuous and positive, and moreover that, for any $y\in S^1$, we have that
$\sum_{x\,\,\text{such that}\,\,\, T(x)= y} \, J(x)=1.$

We say that $\rho$ is a Gibbs probability  for the H\"older continuous potential $\log J$, if
$\mathcal{L}_{\log J}^* \,(\rho)=\rho ,$ where $\mathcal{L}_{\log J}$ is the Ruelle operator for  $\log J$. We call $J$ the Jacobian of $\rho$.

Suppose $\nu=\mu_1*\mu_2$ is the convolution of two Gibbs probabilities  $\mu_1$ and $\mu_2$ associated, respectively, to $\log J_1$ and $\log J_2$. We show that $\nu$ is also Gibbs and its Jacobian $\tilde{J}$ is given by $\tilde{J}(u) = \int J_1(u-x) d \mu_2(x)$

In this case, the entropy $h(\nu)$ is given by the expression
$$ h(\nu) = - \int\,[\,\,\int\, \log \,(\,\int J_1(r+s-x) d \mu_2(x)\,) \, d  \mu_2(r)\,\, ]\,\,d \mu_1 (s).$$
For a fixed $\mu_2$ we consider differentiable variations $\mu_1^t$, $t \in (-\epsilon,\epsilon)$,  of $\mu_1$ on the Banach manifold of Gibbs probabilities, where $\mu_1^0=\mu_1$, and we estimate the derivative of the entropy $h(\mu_1^t * \mu_2)$ at $t=0$.

We also present an expression for the Jacobian of the convolution of a Gibbs probability $\rho$ with the invariant probability with support on a periodic orbit of period two. This expression is based on the Jacobian of $\rho$ and two Radon-Nidodym derivatives.

\end{abstract}

\section{Introduction}

Consider the $2\, x $ (mod 1) transformation $T$ on the unitary circle $S^1$.

All expressions of the form $x+y$ below are consider (mod 1).

Given two probabilities  $\eta $ and $\mu$  on $S^1$ the convolution $ \nu\,=\, \eta * \mu$ is the probability such that for any Borel set $A$ we have
$$ \nu(A)=(\eta * \mu)(A)= \int\,  \mu(A-x) \, d \eta(x).$$

This is the same as saying that
for any continuous function $\phi$
$$\int \phi (z) d \nu(z)=  \int\, (\int \, \phi(y+x)\, d \mu(y))\,d\,\eta(x).$$

Note that if $\mu$ is the Lebesgue probability on $S^1$ then, for any $\nu$ we get $\mu * \nu = \mu$ (just change coordinates).

On the other hand if $\mu =\delta_0$, then, for any $\nu$ we have that $\mu * \nu = \nu$.

\bigskip

Suppose $\mu$ and $\eta$ are $T$ invariant.

Note that for any continuous $\phi$
$$ \int (\phi \circ T)(z)\, \ d \nu (z) =\int\, (\int \, (\phi \circ T )(x +y)\, d \mu(y))\,d\,\eta(x)=\,\,\, $$
$$\int\, (\int \, (\phi (\,T(x) +T(y)\,)\, d \mu(y))\,d\,\eta(x)=$$
$$\int\, (\int \, (\phi (\,T(x) + y\,)\, d \mu(y))\,d\,\eta(x)=$$
$$\int\, (\int \, (\phi (\,x + y\,)\, d \mu(y))\,d\,\eta(x)=\int \phi (z) d \nu(z).$$

Then, it follows that $\nu$ is also $T$-invariant.

The convolution operation is commutative.

A very important contribution to the topic of convolution of invariant probabilities on the circle is \cite{Li}.
By means of combinatorial techniques it was proved there,
among other things, the convergence of the n-convolution of positive entropy measures for
the Lebesgue measure. They also show  that convolution does not decrease entropy (increase in most of the examples). The proofs of our results have a different nature (and they are for a particular family of probabilities).

Related results concerning convolution are \cite{HS}, \cite{Meiri1}, \cite{Meiri2} and \cite{Meiri3}.

\medskip

\begin{definition} \label{saco}
The Jacobian of an invariant measure $\rho$ is the measurable transformation $J_\rho$ such that
$$ \rho(T(A)) = \int_A J_\rho^{\,-1} d \rho,$$
for any Borel set $A$  such that $T|_A$ is injective (see the paragraphs before Proposition 3.4 in \cite{Via}).

\end{definition}

The role of $J_\rho^{\,-1}$ is to provide a formula for the change of variables on the inverse branches of $T$.

A more elegant form of expressing this property for  $J_\rho^{\,-1} $  (in the particular case which is the main interest of our paper ) is via the Ruelle operator. We begin with the "Jacobian" and a posteriori we get the probability. We point out that the Jacobian in \cite{Via} is the inverse of what we call Jacobian here.

\medskip

We assume from now on that $J: S^1 \to \mathbb{R}$ is at least continuous and positive, and such that, for any $y$ we have that
$$\sum_{x\,\,\text{such that}\,\,\, T(x)= y} \, J(x)=1.$$

Given $J$ as above the Ruelle operator $\mathcal{L}_{\log J}$ acts on continuous functions $\varphi$ on the following way:
$\mathcal{L}_{\log J}(\varphi)= \phi$, where
$$\phi(y) =\sum_{x\,\,\text{such that}\,\,\, T(x)= y} \, J(x)\, \varphi(x).$$

The dual $\mathcal{L}_{\log J}^*$ of $\mathcal{L}_{\log J}$ acts on probabilities.

\begin{definition}
We say that $\rho$ is a Gibbs probability (or, a $g$-measure, where $g=\log J$) for the continuous function $J$ if
$$\mathcal{L}_{\log J}^* \,(\rho)=\rho .$$
\end{definition}

The entropy of $\rho$ is given by the Rokhlin formula:
$- \int \log J d \rho$ (see for instance section 9.7 in \cite{VO}). The probability
$\rho$ is the equilibrium probability (maximize pressure) for the potential $\log J$ (see Proposition 3.4 in \cite{PP}).

The Jacobian $J_\rho$ of $\rho$ according to definition \ref{saco} agrees with the above $J$.

In this way is natural to call $J$ the Jacobian of $\rho$.

As an example we mention that for the transformation $T(x)$= $2 x $ (mod $1$) the Lebesgue probability has Jacobian $J$ constant equal to $1/2$.

If $J$ is just continuous it is possible that exists more than one fixed point probability for $\mathcal{L}_{\log J}^*$ (see \cite{BK} and \cite{Quas}). If $J$ is H\"older the fixed point probability is unique.

General references for Jacobians and Thermodynamics Formalism are  \cite{VO}, \cite{MU}, \cite{PP}, \cite{Par} and \cite{Rue}. We use the dynamics of the doubling map on an essential way. The possible extension to expanding transformations on the circle would require a good meaning for translation on the circle which is at the same time compatible with
the distance among preimages of a general point.

 In the  section \ref{Gi} we will consider convolution of two Gibbs probabilities.  We estimate the entropy of the convolution of two Gibbs probabilities (see Theorem \ref{est1}). We also show for the case of Gibbs probabilities that, if $ \nu=\mu_1 * \mu_2$, then,  $h(\nu) \geq h(\mu_2)$ (see Theorem \ref{est3}). This result appears in a more general setting in \cite{Li}. We do not use here the Hausdorff dimension as a tool in our proof.

We will present in section \ref{per}  an explicit expression for the Jacobian of the probability obtained by the
convolution of a Gibbs probability and a periodic orbit of period two
 (see expression (\ref{cio1})).

We also show examples of Gibbs  probabilities $\mu$ where the convolution of $\mu$ with a periodic orbit of period two results in the same probability $\mu$ (see the class of potentials $\mathcal{S}$ defined by expression (\ref{haha})).

In section \ref{dif}  we analyze the following problem:
for a fixed $\mu_2$ consider differentiable variations $\mu_1^t$, $t \in (-\epsilon,\epsilon)$,  of $\mu_1$ on the Banach manifold of Gibbs probabilities, where $\mu_1^0=\mu_1$. How can one   estimate the derivative of the entropy $h(\mu_1^t * \mu_2)$ at $t=0$? On this direction see Proposition \ref{ert}.

In the appendix we consider the following problem:
suppose   $J_1$ and $J_2$  are the H\"older Jacobians and they are such that: $J_2\geq J_1$, when $J_1\geq1/2$, and $J_2\leq J_1$, when $ J_1\leq 1/2$.
Denote $\mu_i$ the Gibbs probability associated to the potential $\log J_i$, $i=1,2$. We show that $h(\mu_1)\geq h(\mu_2)$ (see Proposition \ref{por}). This problem is related to questions raised in section \ref{per}.

 The PhD thesis \cite{Uggi} and \cite{Bar} consider several  properties for  the convolution of invariant probabilities for the symbolic space setting. An appropriate structure have to be considered for replacing the sum translation on the circle. These  works do not consider results similar to ours.

We thanks L. Cioletti, P. Giulietti and B. Uggioni for helpful conversations on the topic of convolution of  invariant probabilities.

 \section{Convolution of Gibbs probabilities} \label{Gi}

Suppose $J_2$ is a  H\"older Jacobian and $J_1$ is a Jacobian which is just continuous. As we said $J_i: S^1 \to \mathbb{R}$, $i=1,2$, are such that $\mathcal{L}_{\log J_i}^* ( \mu_i)=\mu_i$. The probability $\mu_2$ is invariant, ergodic and has support on $S^1$.

We want to estimate analytical properties of the probability
$\nu=\mu_1 * \mu_2$.

A natural question is to ask if there exists an explicit expression for the
Jacobian $\tilde{J}$, such that,
$$\mathcal{L}_{\log \tilde{J}}^* \,(\nu)=\nu $$ in terms of $J_1,J_2$.

\begin{theorem} \label{est1}
Suppose $J_2$ is a  H\"older  and $J_1$ is continuous.
Then,  the Jacobian  $\tilde{J}$  of $\nu=\mu_1 * \mu_2$ satisfies for any $u$ the expression
\begin{equation} \label{ot0} \tilde{J}(u) = \int J_1(u-x) d \mu_2(x)
\end{equation}

and, therefore
$$ h(\nu) = - \int \log \tilde{J} (u) d \nu(u) = - \int\, \log \,(\,\int J_1(u-x) d \mu_2(x)\,) \, d  \nu(u)=$$
\begin{equation} \label{ot1}- \int\,[\,\,\int\, \log \,(\,\int J_1(r+s-x) d \mu_2(x)\,) \, d  \mu_2(r)\,\, ]\,\,d \mu_1 (s). \end{equation}

\bigskip

\end{theorem}

\medskip

In the proof of this theorem we just  need
to use the fact that $\mathcal{L}_{\log J_1}^* \,(\mu_1)=\mu_1 $ and it is not required that $\mu_1$ is the limit of the $\rho_n$, $n \in \mathbb{N}$,  defined by (\ref{hu}). However, this property is required for $\mu_2$. The proof will be done later.

Note that by the commutativity of the convolution we get that the above defined function $\tilde{J}$ is Holder if either $J_1$ or $J_2$ is Holder.

\medskip

\begin{corollary}  Suppose $\mu_1$ has a Jacobian $J_1$ which is continuous and $\mu_2$ is any invariant probability.
Then,  the Jacobian $\tilde{J}$ of $\nu=\mu_1 * \mu_2$ satisfies for any $u\in S^1$ the expression
\begin{equation} \label{ot7} \tilde{J}(u) = \int J_1(u-x) d \mu_2(x)
\end{equation}

\end{corollary}
{\bf Proof:}  Any invariant probability $\mu_2$ can be weakly approximated by Gibbs states $\mu_2^n$, $n \in \mathbb{N}$ (see for instance Theorem 8 page 536 in \cite{L3}).

The function $\rho \to \mu_1 * \rho $ is continuous in the weak topology.

Then, the Jacobian $\tilde{J}_n$ of $\nu_n=\mu_1 * \mu_2^n$ converges to the function $\tilde{J}(u)= \int J_1(u-x) d \mu_2(x).$ Indeed,  $x \to J_1(u-x)$ is a continuous function depending $C^0$ on $u$.

The function $\tilde{J}$ is continuous positive and satisfies $\tilde{J}(x_1) + \tilde{J}(x_2)=1$, if $T(x_1)=T(x_2).$

In order to show that $\tilde{J}$ is the Jacobian of $\nu=\mu_1 * \mu_2$ consider any arbitrary continuous function $\varphi$.

Then,
$$ \int \mathcal{ L}_{ \log \tilde{J} } \, (\varphi)\,(z) d \nu(z) = \int  \sum_{T(w)=z}  \,[\,\int J_1(w-x) d \mu_2(x) \,] \,\varphi(w) \, d \nu(z)=$$
$$  \int  \sum_{T(w)=z}  \,[\,\int J_1(w-x) d \mu_2(x) \,] \,\varphi(w) \, d (\mu_1 * \mu_2)(z)=$$
$$ \lim_{n \to \infty} \int \,\sum_{T(w)=z}  \,[\,\int J_1(w-x) d \mu_2^n(x) \,] \,\varphi(w) \, d (\mu_1 * \mu_2^n)(z)=$$
$$ \lim_{n \to \infty} \int \,  \,\varphi(z) \, d (\mu_1 * \mu_2^n)(z)=  \int \,  \,\varphi(z) \, d (\mu_1 * \mu_2)(z)= \int \,  \,\varphi(z) \, d \nu(z).$$

Therefore, $\mathcal{ L}_{ \log \tilde{J} }^* ( \nu)= \nu.$ Finally, from Theorem \ref{est1} we get that $u \to \int J_1(u-x) d \mu_2(x)$ is the continuous Jacobian of $\nu=\mu_1 * \mu_2.$

\qed

\medskip

\begin{corollary} \label{est2} Suppose $J_1$ is a  H\"older, $J_2$ is continuous and $\mu_2$ is the limit of the probabilities $\rho_n$ defined on (\ref{hu}).
Then, the Jacobian
$\tilde{J} $ of $\mu_1 * \mu_2$
is  H\"older and has the same H\"older constant. This means that convolution regularizes Jacobian.

\end{corollary}
{\bf Proof:}

As we mention in the remark at the end of this section the expression $\tilde{J}(u) = \int J_1(u\,-x) d \mu_2(x)$ is true.

Suppose $0<\alpha\leq 1$ and $K$ are such that for any $r,s$ we have
$$ |J_1(r)-J_1(s)| \leq K\, |r-s|^\alpha,
$$
then, for any $u_1,u_2$
$$ |\,\int J_1(u_1 - x)\, d\mu_2(x) \,-\,\int J_1(u_2 - x)\, d\mu_2(x)\,| \leq \, $$
$$\,\int \,|J_1(u_1 - x)\, \,-\,J_1(u_2 - x)\,|\, d\mu_2(x)\leq     K\, |u_1-u_2|^\alpha.
$$

\qed

\medskip

It is known from Lemma 9.2 (or, Corollary 9.3) in \cite{Li} that convolution increase entropy, that is, $h(\mu_1 * \mu_2) \geq h(\mu_2)$. The proof in \cite{Li} basically
use the fact that $HD(\mu)= \frac{ h(\mu) }{\log 2}$ and simple properties of the Hausdorff dimension of an invariant probability. We will present a direct proof without using Hausdorff dimension for the case of Gibbs probabilities.
We point out that Gibbs probabilities are dense in the set of invariant probabilities (see for instance Theorem 8 page 536 in \cite{L3}).

\begin{theorem} \label{est3} Suppose $J_1$ and $J_2$ are  H\"older Jacobians. Denote by $\mu_1$ and $\mu_2$ the corresponding  Gibbs probabilities. If $ \nu=\mu_1 * \mu_2$, then,  $h(\nu) \geq h(\mu_2)$. Moreover, we have that $h(\nu) > h(\mu_2)$, unless $\mu_1$ or $\mu_2$ is the Lebesgue probability.

\end{theorem}

{\bf Proof:}
It is known from \cite{Lop} (or, \cite{LMMS} for a more general statement) that when $\mu_2$ has a H\"older Jacobian we get
$$ h(\mu_2) =  \inf_{v >0,\,\,v\,\text{H\"older}} \,\int \log \bigg(\frac{\mathcal{L}_{0} v(s)}{v(s)}\bigg) d\mu_2(s) .$$
where for any $s$ we have $ \mathcal{L}_{0} v(s) = v(s_1) + v(s_2).$ This condition can be relaxed assuming that $v$ is just continuous (indeed, one can check  that the proof of Lemma 2 in \cite{LMMS} applies to continuous potentials).

We will show that there exists $u$ such that

$$h(\nu) \geq \int \log \bigg(\frac{\mathcal{L}_{0} u(s)}{u(s)}\bigg) d \mu_2(s)= $$
$$\int \log (\mathcal{L}_{0} u(s)) d \mu_2(s)\,-\,\int \log  u(s) d \mu_2(s) .$$

More precisely we will exhibit a H\"older continuous function  $u$ such that
$$ -\int \log u (s)\, d \mu_2(s) = h(\nu),$$
and,
moreover that
$$\int \log (\mathcal{L}_{0} u(s))\, d \mu_2(s)\leq 0.$$

From (\ref{ot1}) we have that
$$ h(\nu) =    - \int\,[\,\,\int\, \log \,(\,\int J_1(r+s-x) d \mu_2(x)\,) \, d  \mu_2(r)\,\, ]\,\,d \mu_1 (s)\,= $$
$$- \int\,[\,\,\int\, \log \,(\,\int J_1(r+s-x) d \mu_2(x)\,) \, d  \mu_1(r)\,\, ]\,\,d \mu_2 (s). $$
Then, taking
$$u(s) = e^{ \,\,\int\, \log \,(\,\int J_1(r+s-x) d \mu_2(x)\,) \, d  \mu_1(r)\,\,  },$$
we just have to show  that $\mathcal{L}_{0} (u(s))\leq 1.$

Suppose $s_1$ and $s_2$ are the two preimages of $s$, then,
$$\int \int J_1(r+s_1-x)\,d \mu_2(x) d  \mu_1(r) + \int\int J_1(r+s_2-x)\,d \mu_2(x)\,  d  \mu_1(r)=1.$$

From Jensen inequality we get that
$$u(s_1) + u(s_2) =$$
$$ e^{ \,\,\int\, \log \,(\,\int J_1(r+s_1-x) d \mu_2(x)\,) \, d  \mu_1(r)\,\,  }+ e^{ \,\,\int\, \log \,(\,\int J_1(r+s_2-x) d \mu_2(x)\,) \, d  \mu_1(r)\,\,  }\leq$$
$$ e^{ \,\log [\,\int\, \,(\,\int J_1(r+s_1-x) d \mu_2(x)\,) \, ]\,d  \mu_1(r)  }+ e^{ \,\, \log \,[\,\int\,(\,\int J_1(r+s_2-x) d \mu_2(x)\,) \,]\, d  \mu_1(r)\,\,  }=$$
$$\int\int J_1(r+s_1-x)\,d \mu_2(x) d  \mu_1(r) + \int \int J_1(r+s_2-x)\,d \mu_2(x)\, d  \mu_1(r)=1.$$

Note that if for some $s$ we have that $\mathcal{L}_{0} u(s)< 1$, then, as $\mu_2$ has full support (see \cite{PP}), we have strict inequality $h(\nu) >h(\mu_2)$.
In order to prevent this from happening it is required that for any $s$
$$\int\, \log \,(\,\int J_1(r+s-x) d \mu_2(x)\,) \, d  \mu_1(r)= $$
$$\log \,\int\,(\,\int J_1(r+s-x) d \mu_2(x)\,) \, d  \mu_1(r).$$

Note that when $J=1/2$ (the Lebesgue probability) then the above equality is true.

On the other hand, if the above equality is true for any $s$ then $J_1$ is constant (equal to $1/2$).
Indeed, it is know that the Jensen inequality is an equality just when all weights are equal.
It follows that $\int J_1(r+s-x) d \mu_2(x)\,$ is constant independent of $r$ and $s$. As $\mu_2$ has full support we get that $J_1$ is constant.

\qed

\medskip
We will show later in section \ref{per} that  there are examples in which the convolution of a Gibbs probability with a probability with support on a periodic orbit results on the initial Gibbs probability.

\begin{theorem} \label{est4} Suppose $\mu$ is Gibbs probability for a H\"older Jacobian $J$. For each $n \in \mathbb{N}$ denote
$ \nu_n=\underbrace{\mu * \mu *....* \mu}_n$, then,  $\lim_{n \to \infty} \nu_n$   is the Lebesgue probability\,\,\,

\end{theorem}

{\bf Proof:} If $\mu$ is the Lebesgue probability there is nothing to prove.

 The sequence of probabilities $\nu_n$, $n \in \mathbb{N}$,  has a convergent subsequence, $\nu_{n_k}$, $k \in \mathbb{N}$. Suppose $\lim_{k \to \infty} \nu_{n_k}= \rho$ and $\rho$ is not Lebesgue probability.

 Denote by $\overline{J}_k $ the Jacobian of $\nu_{n_k}.$ The sequence $\overline{J}_k $, $k \in \mathbb{N}$, is equicontinuous and bounded
 by Theorem \ref{est2}. Then, by Arzela-Ascoli theorem there exist an uniform limit $\overline{J}_\infty $  (which is H\"older) of a subsequence of $\overline{J}_k $, $k \in \mathbb{N}$.

\medskip
\begin{remark} \label{que}. By weak$*$ topology one can show that the Jacobian of such probability  $\rho$ is exactly $ \overline{J}_\infty $ \,\,\label{que}.
\end{remark}
 \medskip

  Denote by $\alpha$  the supremum of the entropy of $h(\rho)$ among the possible $\rho$ obtained by convergent subsequences, $\nu_{n_k}$, $k \in \mathbb{N}$.

  We claim that one $\hat{\rho}$ of such possible $\rho$ attains the supremum.

  Consider a sequence of $\hat{\rho_r}$, $r \in \mathbb{N}$ of such possible limit of subsequences $\nu_{n_k}^r$, $r \in \mathbb{N}$, $n \in \mathbb{N}$ such that
 $$\lim_{r \to \infty} \hat{\rho_r}= \overline{\rho}, $$
 and
  $$\lim_{r \to \infty} h(\hat{\rho_r})= \alpha. $$

Then, it is possible to get a sequence $\nu_{n_{k(r)}}^r$ such that
$$\lim_{r \to \infty} \nu_{n_{k(r)}}^r= \overline{\rho}, $$
and
$$\lim_{r \to \infty} h(\nu_{n_{k(r)}}^r)= \alpha. $$

As the entropy is lower semicontinuous we get that $h(\overline{\rho})=\alpha$. By  Remark \ref{que}  we get that
$\overline{\rho}$ has a H\"older Jacobian.

  Suppose $\alpha<\log 2$. Then,  we get by  Theorem  \ref{est3} that $\mu * \overline{\rho}$ has  bigger entropy than $\overline{\rho}$. If $\lim_{k \to \infty} \nu_{n_k}=  \hat{\rho}$ then
  $\lim_{k \to \infty} \mu *\nu_{n_k}=  \mu *\hat{\rho}$ and this is a contradiction.

This proves that $\alpha=\log 2$ and, by monotonicity of the entropy function along  the sequence $\nu_n$, that the unique maximal entropy measure (Lebesgue) is the weak limit of $\nu_n$.

\qed

\bigskip

In \cite{Li} the authors proved convergence to Lebesgue measure for concatenations $\mu_n*...*\mu_2*\mu_1$  (invariant measures) with some bound on their entropy. In the moment we don't know  how to get this kind of result with our methods.

\bigskip

Suppose $J$ is the H\"older Jacobian of the probability $\rho.$

Consider for each $n\in \mathbb{N}$, the probability
\begin{equation} \label{hu} \rho_n= \sum_{j=1}^{2^n}  \delta_{ \frac{j}{2^n}} \Pi_{k=0}^{n-1}  J (\, T^k (\, \frac{j}{2^n}\, )\,)  \end{equation}
which is not $T$ invariant.

Note that $T^n ( \frac{j}{2^n})=0$.

 If $\lim_{n \to \infty} \rho_n = \rho$, then $\rho$ will satisfy the equation $\mathcal{L}_{\log J}^* \,(\rho)=\rho.$ (see Ruelle Theorem \cite{PP})

\begin{remark} \label{hu1} Note  that $\rho_n =(\mathcal{L}_{\log J}^*)^n (\delta_0)$.  In the case $J$ is Holder it is a classical result
that $\rho_n \to \rho$  (see Ruelle Theorem 2.2 (iv) in \cite{PP}). A more strong claim is Theorem 1.1 in \cite{KLS} which shows convergence on the $1$-Wassertein distance.
\end{remark}

 In the case $J$ is continuous  we will assume here that such limit exists $\rho_n \to \rho$ and we point out that this limit is a Gibbs state for $J$.

 If $J$ is H\"older such limit exist and it is the only fixed point of $\mathcal{L}_{\log J}^*$.

\bigskip

Now we will begin the proof of Theorem \ref{est1}.
From now on we denote $\mu_1=\mu$ and $J_1=J$.
We want to determine $\tilde{J}$ from $J$.

Denote $J_2^n (j):=  \Pi_{k=0}^{n-1}  J_2 (\, T^k (\, \frac{j}{2^n}\, ))$, $j=1,2,..,2^n$.

Now we will consider $\rho_n$ when $J=J_2$.
That is, $\rho_n$ is the probability
$$ \rho_n= \sum_{j=1}^{2^n}  \delta_{ \frac{j}{2^n}}   J_2^n (j)  ,$$
which is not $T$ invariant.

It is known (see Remark \ref{hu1}) that
$$ \lim_{n \to \infty} \rho_n = \mu_2$$

\medskip

It is natural to consider in our reasoning
the convolution $\mu * \rho_n= \nu_n$, $n \in \mathbb{N}$, because $\nu_n \to \nu$, when $n \to \infty.$ We denote by $\tilde{J}_n$ the Jacobian of the (in principle) non invariant probability $\nu_n$.

\medskip

Suppose $ y$ is such that $ \frac{k}{2^n } \leq y< \frac{k+1}{2^n} .$ For fixed $j$, what is the range of $x$ such that $ y = \frac{x}{2} + \frac{j}{2^n}.$ The answer is $ \frac{k-j}{2^{n-1}} \leq x< \frac{k-j+1}{2^{n-1}}.$

\medskip

We will show later that for $u\in [\frac{v}{2^{n+1}}, \frac{v+1}{2^{n+1}}) $
$$ \tilde{J}_n(u)=  \sum_{ j\,=1}^{2^{n-1}}\, \,  J_2^n (j) \,J( u- \frac{j}{2^n})    . $$
\medskip

Note that for a continuous function $f$ we get

$$\int f(z) d (\mu * \rho_n) (z) =
 \int \,( \int f(x+y) d \mu (x)\, )\, d \rho_n (y)\, ) =$$
$$ \int\,( \int \mathcal{L}_{\log J} (f(x+y)) d \mu (x)\,) d \rho_n (y) =$$
$$ \sum_{ j=1}^{2^n}  \Pi_{w=0}^{n-1}  J_2 ( T^w (\, \frac{j}{2^n})) \,\int [ J(x/2)  (f(x/2\,+ \frac{j}{2^n} )) + J((x+1)/2)  (f(\,(x+1)/2+ \frac{j}{2^n}) )]  d \mu (x))= $$
 \begin{equation}  \label{w1}\sum_{ j=1}^{2^n} \Pi_{w=0}^{n-1}  J_2 ( T^w ( \frac{j}{2^n}))  \int [ J(x/2)  f(x/2+ \frac{j}{2^n} ) + J(x/2+1/2)  f(x/2 +1/2+ \frac{j}{2^n})]  d \mu (x) )= \end{equation}
  \begin{equation}  \label{w2}\sum_{ j=1}^{2^n} \Pi_{w=0}^{n-1}  J_2 ( T^w ( \frac{j}{2^n} )) \int [  \tilde{J}_n(x/2+\frac{j}{2^{n+1}})  f(x/2 +\frac{j}{2^{n+1}}) + \tilde{J}_n(x/2+\frac{j}{2^{n+1}}+1/2)  f(x/2+\frac{j}{2^{n+1}}+1/2)]  d \mu (x) = \end{equation}
$$  \sum_{ j=1}^{2^n}\Pi_{w=0}^{n-1}  J_2 (\, T^w (\, \frac{j}{2^n}\, )) \int [ \tilde{J}_n((x+\frac{j}{2^n})/2)  f((x+\frac{j}{2^n})/2) + \tilde{J}_n((x+\frac{j}{2^n}+1)/2)  f(\,(x+\frac{j}{2^n}+1)/2\,)\, ]  d \mu (x) = $$
$$  \int \mathcal{L}_{\log \tilde{J}_n} (f) (x+y)  d \mu(x)  d \rho_n (y)=  $$
$$ \int \mathcal{L}_{\log \tilde{J}_n} (f) (z)  d (\mu * \rho_n ) (z)= \int f (z)  d (\mu * \rho_n ) (z).  $$

\medskip

Below we consider any $j$ modulo $2^n$.

\medskip

Suppose $f$ is a function with support on $[\frac{v}{2^{n+1}}, \frac{v+1}{2^{n+1}})$, $0\leq v \leq 2^{n+1}-1$.

In  figure 1 we consider the case $n=2$, and one can see, for instance, on the interval $[\frac{3}{2^3},\frac{4}{2^3})=[\frac{v}{2^{n+1}}, \frac{v+1}{2^{n+1}})$,  that { two branches $ x/2$ and $x/2+ \frac{1}{2^2}$, have projections over  $(\frac{3}{2^3},\frac{4}{2^3})$, (using the red color - this corresponds $j=1,3$ and to left hand side of (\ref{w1}) \,)}  and, moreover, $ x/2 $, $x/2 +  \frac{1}{2^3}$,  $x/2+  \frac{2}{2^3}$ $x/2 +  \frac{3}{2^3}$ (using the red and the blue color - this corresponds to $j=1,2,3,4$ and  (\ref{w2})\,), also have projections over  $(\frac{3}{2^3},\frac{4}{2^3})$.

{ In the general case for the interval  $[\frac{v}{2^{n+1}}, \frac{v+1}{2^{n+1}})$,  $ v=0,1,..., 2^{n+1}-1$, } we have to consider for (\ref{w1})

a) for $v$= even it is required a range of values $j$ where $j=\frac{v- t \,2}{2}$, $t=0,..,2^{n-1}-1$  (for the left hand side of (\ref{w1}) ). Moreover, for the right hand side of  (\ref{w1}) we will need the values of $j= \frac{v-2t}{2}- 2^{n-1}$.

b) for $v$=odd it is required a range of values
$j$ where $j=\frac{v-1- t \,2}{2}$, $t=0,..,2^{n-1}-1$  (for the left hand side of (\ref{w1}) ). Similar as above for the right hand side.

{This means the total of $2^{n}$ possible values of $\frac{j}{2^n}$ in each case a) or b). We use this identification of $t$ and $j$ on future expressions.}

\medskip

{
For the interval  $[\frac{v}{2^{n+1}}, \frac{v+1}{2^{n+1}})$  we have to consider at same time the both expressions (left and right) of the sum for (\ref{w2}). Note that $v$ ranges on $0,1,..,2^{n+1}$. Given $j$, there exists a $j_0\in\{1,2,..,2^n\}$ such that either $j+j_0=v$ or $j+j_0= v-2^n$. Each $j\in\{1,2,..,2^n\}$ can not satisfy both conditions at same time. Any $j\in\{1,2,..,2^n\}$ will satisfy one of the conditions. In this way all $j$ will be used when considering together the left and right side of (\ref{w2}).}


\begin{figure}
\begin{center}
\includegraphics[width=4.5in]{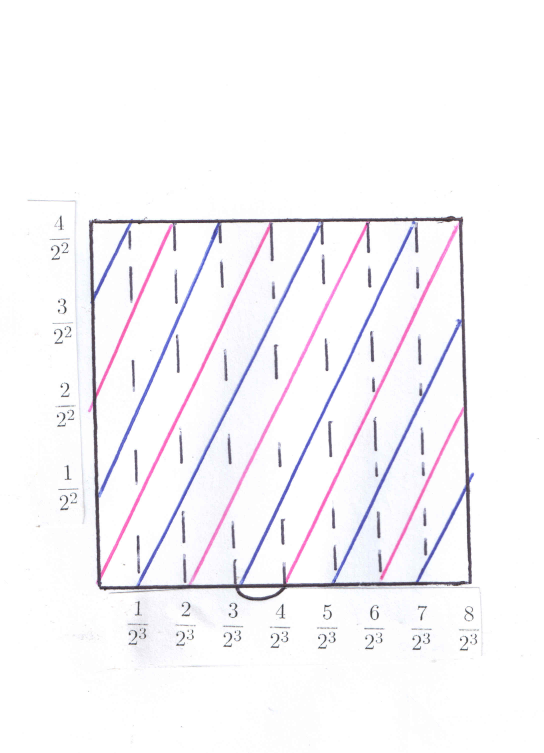}
\caption{ The case $n=2$}
\end{center}
\end{figure}



\medskip

We assume now that $ v=0,1,..., 2^{n+1}-1$ is even.

In this case we consider the two terms of (\ref{w2}):

\begin{equation}  \label{aa1}
\textcolor{black}{\sum_{j\,\,\text{such that} \,j+j_0= v\, \text{for some}\, j_0} \Pi_{w=0}^{n-1}  J_2 (\, T^w (\, \frac{j}{2^n}\, )\,) \int_{ \frac{v-j}{2^n}}^{  \frac{v-j+1}{2^n}} \,\tilde{J}_n (\frac{x}{2}\,+\,\frac{j}{2^{n+1}}  )\, f(\frac{x}{2}\,+\,\frac{j}{2^{n+1}}  )\, d \mu(x)\,}+\,\end{equation}
\begin{equation}  \label{aa2}
\textcolor{black}{\sum_{j\,\,\text{such that} \,j+j_0= v-2^n\, \text{for some}\, j_0 }  \Pi_{w=0}^{n-1}  J_2 ( T^w ( \frac{j}{2^n} ))  \int_{ \frac{v-j-2^n}{2^n}}^{  \frac{v-j-2^n +1}{2^n}} \,\tilde{J}_n (\frac{x}{2}+\frac{j+ 2^n }{2^{n+1}}  )\,f(\frac{x}{2}+\frac{j + 2^n}{2^{n+1}}  ) \,d \mu(x)}=\end{equation}
\begin{equation}  \label{aa3}
\textcolor{black}{\sum_{j\,\,\text{such that} \,j+j_0= v\, \text{for some}\, j_0 }\Pi_{w=0}^{n-1}  J_2 (\, T^w (\, \frac{j}{2^n}\, )\,) \int_{ \frac{v}{2^{n+1}}}^{  \frac{v+1}{2^{n+1}}} \,\frac{\tilde{J}_n (u)}{J_n(u)}\, f(u)\, d \mu(u)\,}+\,\end{equation}
\begin{equation}  \label{aa4}
\textcolor{black}{\sum_{j\,\,\text{such that} \,j+j_0= v-2^n\, \text{for some}\, j_0  } \Pi_{w=0}^{n-1}  J_2 (\, T^w (\, \frac{j}{2^n}\, )\,) \int_{ \frac{v}{2^{n+1}}}^{  \frac{v+1}{2^{n+1}}} \,\frac{\tilde{J}_n (u)}{J_n(u)}\, f(u)\, d \mu(u)}=\end{equation}
\begin{equation}  \label{aa5}
{      \,    \int_{ \frac{v}{2^{n+1}}}^{  \frac{v+1}{2^{n+1}}} \,\frac{\tilde{J}_n (u)}{J_n(u)}\, f(u)\, d \mu(u)}.\end{equation}
\bigskip

Assume that $ v=0,1,..., 2^{n+1}-1$ is even.  In this case we consider the two terms of (\ref{w1}):

$$ \sum_{ 0\leq t\leq 2^{n-1}-1\,}\,  \Pi_{w=0}^{n-1}  J_2 (\, T^w (\, \frac{\frac{v-2t}{2}}{2^n}\, )\,) \, {\int_{\frac{2\, t}{2^n}}^{\frac{2\,t +1}{2^n}}  J(x/2)  f(x/2\,+ \frac{\frac{v-2t}{2}}{2^n} )  d \mu (x) \,}+$$

$$  \sum_{ 0\leq t\leq 2^{n-1}-1\,}\,  \Pi_{w=0}^{n-1}  J_2 (\, T^w (\, \frac{\frac{v-2t}{2}- 2^{n-1}}{2^n}\, )\,)\, {\int_{\frac{2\, t}{2^n}}^{\frac{2\,t +1}{2^n}}  J(x/2\,+1/2)  f(\,x/2\, +1/2\,+\, \frac{\frac{v-2t}{2}- 2^{n-1}}{2^n}) d \mu (x)}= $$
$$ \sum_{ j\,\,=1}^{2^{n-1}}\, \,  \Pi_{w=0}^{n-1}  J_2 (\, T^w (\, \frac{j}{2^n}\, )\,) \, {\int_{\frac{v}{2^{n+1}}}^{\frac{v+1}{2^{n+1} }} \frac{J( u- \frac{j}{2^n})}{J_n(u)}  f(u) \,d \mu (u)}  .$$

Therefore, when $v$ even we get that for $u\in [\frac{v}{2^{n+1}}, \frac{v+1}{2^{n+1}}) $
$$ \tilde{J}_n(u)=  \sum_{ j\,=1}^{2^{n-1}}\, \,  \Pi_{w=0}^{n-1}  J_2 (\, T^w (\, \frac{j}{2^n}\, )\,) \, \textcolor{black}{ J( u- \frac{j}{2^n})  }   $$
\medskip
A similar result is true when $v$ is odd.

\bigskip

Remember that $ \lim_{n \to \infty} \rho_n = \mu_2$.

As $J_2$ is a  H\"older function, given any $x_0\in S^1$ we have that $\lim_{n \to \infty} \mathcal{L}^n_{\log J_2} (g) \,(x_0) = \int g(z) d \mu_2(z)$. Then, we consider for $u\in S^1$ fixed, the function $x \to g(x)= J( u -x) .$

Then, taking $x_0=0$ we get
$$ \tilde{J} (u)=\lim_{n \to \infty} \, \tilde{J}_n (u) =   \lim_{n \to \infty} \,  \mathcal{L}^n_{\log J_2} (g) \,(0) =\int J(u-x) d \mu_2(x).$$

In this way

$$ h(\nu) = - \int \log \tilde{J} (u) d \nu(u) = - \int\, \log \,(\,\int J(u-x) d \mu_2(x)\,) \, d  \nu(x). $$

\bigskip

\begin{remark} \label{why} Note also that if $J_2$ is continuous and satisfies the hypothesis of being the limit of $\rho_n$, $n \in \mathbb{N}$, the same expression $ \tilde{J} (u)=\int J(u-x) d \mu_2(x)$ obtained above is also true.
\end{remark}

\medskip

In the case $J$ is constant $J=1/2$ we get that $\tilde{J}_n=1/2$. In this way if $\mu$ is the Lebesgue probability, then, $\nu_n=\mu * \rho_n$ is also Lebesgue probability.

Note that for any $u$ we have that $  \tilde{J}_n(u)+  \tilde{J}_n(u+1/2)=1$.  In this way the probability $\nu_n$ is invariant for the $T$. Then, we get that the convolution of any invariant probability $\mu$ with $\rho_n$ (not invariant) is invariant.

The entropy of $\nu_n$ satisfies $h(\nu_n)= - \int \log \tilde{J}_n d \nu(n) $.

\bigskip


\section{Convolution of Gibbs probability and a periodic orbit of period two} \label{per}

In this section we consider the convolution of a Gibbs probability with a probability with support on an orbit of period two.

Suppose the Jacobian $J: S^1 \to \mathbb{R}$ is such that $\mathcal{L}_{\log J}^* ( \mu)=\mu$.

Consider now $\rho= \frac{1}{2} ( \delta_{1/3} + \delta_{2/3})$ and we want to analyze properties of
$\nu=\mu * \rho$.

We denote the Jacobian of $\nu$ (in the sense of Definition \ref{saco})  by  $\tilde{J}$.
We have to understand  in this case the corresponding change of coordinates on the inverse branches.

In other words, we want  to express the $\tilde{J}$, such that,
\begin{equation} \label{rrr}  \mathcal{L}_{\log \tilde{J}}^* \,(\nu)=\nu
\end{equation}
in terms of $J$, $\rho$ and $\mu$.

We will present an explicit expression for $\tilde{J}$ in terms of $J$ and two more Radon-Nikodym derivatives (see expression (\ref{cio1})). This will provide a formula for the entropy of $\mu * \rho$ (see (\ref{cio3})).

We will also show that there exist Gibbs probabilities $\mu$ satisfying  $\mu=\nu= \mu * \rho.$
Jacobians described by equation (\ref{haha}) satisfy this property. For these examples, of course, the entropy does not increase by convolution.

\medskip

About question (\ref{rrr}) the main property for $\tilde{J}$ is: for any continuous function $f$
$$ \int \mathcal{L}_{\log \tilde{J}} (f) (z)  d (\mu * \rho) (z) = \int f(z) d (\mu * \rho) (z) .  $$

In this way $\mu * \rho$ is a fixed point for $\mathcal{L}_{\log \tilde{J}}^*$.

Remember that when $\tilde{J}$ is H\"older the fixed point probability is unique.

For a continuous function $f$ we get
 $$\int f(z) d (\mu * \rho) (z) = \int f(x+y) d \mu (x)\, d \rho(y)=$$
 $$ 1/2\,( \int f(x+1/3) d \mu (x)\, +\, \int f(x+2/3) d \mu (x)\, ) =$$
$$ 1/2\,( \int \mathcal{L}_{\log J} (f(x+1/3)) d \mu (x)\, +\, \int \mathcal{L}_{\log J}( f(x+2/3)) d \mu (x)\, ) =$$
$$        1/2\,(\, \int [ J(x/2)  (f(x/2\,+\,1/3)) + J((x+1)/2)  (f(\,(x+1)/2\,+\,1/3)\, )]  d \mu (x) + $$
$$ \int [ J(x/2)  (f(x/2\,+\,2/3)) + J((x+1)/2)  (f(\,(x+1)/2\,+\,2/3)\, )]  d \mu (x) \,)= $$
$$        1/2\,(\, \int [{ J(x/2) \, f(x/2\,+\,1/3\,)} + { J((x+1)/2) \, f(\,x/2\,+\,5/6\,)}\, \,]\,  d \mu (x) + $$
$$ \int [ {J(x/2)\,  f(x/2\,+\,2/3\,)} + {J((x+1)/2) \, f(\,x/2\,+\,1/6\,)}\, ] \, d \mu (x) \,). $$
\medskip

On the other hand

$$ \int \mathcal{L}_{\log \tilde{J}} (f) (z)  d (\mu * \rho) (z)=  $$
$$  \int \mathcal{L}_{\log \tilde{J}} (f) (x+y)  d \mu(x) \,  d \rho (y)=  $$
$$ \int [ \tilde{J}((x+y)/2)  (f(\,(x+y)/2)) + \tilde{J}((x+y+1)/2)  (f(\,(x+y+1)/2\,)\, )]  d \mu (x) \rho(y)\,= $$
$$ 1/2\,(\,\int [\, \tilde{J}((x+\,\,1/3)/2)  \,f(\,(x+\,\,1/3)/2)\, + \tilde{J}(\,(x+\,\,1/3\,\,+1)/2)  \, f(\,(x+\,\,1/3\,\,+1)/2\,)\, ]  d \mu (x) \,+ $$
$$ \int [ \,\tilde{J}(\,(\, x\,+\,2/3\,)/2\,)  \,f\,(\,(\,x\,+\,\,2/3\,\, )/2\,) + \tilde{J}(\,(\,x +\,2/3\,+1\,)/2\,)  \,f(\,(\,x\,+\,\,2/3\,+\,1\, )/2\,)] \, d \mu (x) \,\,. $$

The above means that it is required that for any continuous $f$
$$        \, \int [ J(x/2) \, f(x/2\,+\,1/3\,) +  J((x+1)/2) \, f(\,x/2\,+\,5/6\,)\, \,]\,  d \mu (x) + $$
$$ \int [ J(x/2)\,  f(x/2\,+\,2/3\,)+ J((x+1)/2) \, f(\,x/2\,+\,1/6\,)\, ] \, d \mu (x) \,= $$
$$ \,\,\int [ \,\tilde{J}(x/2\,+\,\,1/6\,)  \,f(\,x/2+\,\,1/6) + \,\tilde{J}\,(x/2\,+\,\,2/3\,\,)  \,f(\,x/2\,+\,\,2/3\, )\,]  \,d \mu (x) \,+ $$
 \begin{equation} \label{klm34} \int [ \,\tilde{J}(x/2+\,\,1/3\,) \, f\,(x/2+\,\,1/3\,)\, +  \, \tilde{J}(x/2\,+\,5/6\,) \,f(\,x/2+\,\,5/6\, )\,] \, d \mu (x) \,. \end{equation}

\medskip

\medskip

\subsection{An explicit expression for the convolution in the case of Gibbs probabilities for H\"older Jacobians}

\noindent

Consider a H\"older Jacobian $J$ on $S^1$   and suppose that  $J$ is  the Jacobian  of $\mu$.

 $\tilde{J}$ denotes the Jacobian of $\nu= \mu * \rho$.

We will not be able to show that $\tilde{J} $ is continuous (just measurable). Anyway, we denote $\tilde{J} $ as the Jacobian of $\nu$ in the sense of Remark \ref{saco}.

In this subsection we want to show an explicit expression (in terms of certain Radon-Nidodym derivatives) for  $\tilde{J}$ (see (\ref{cio1})). In order to get that  we will have to use equation (\ref{klm34}).

Denote by $\mu_1$ the probability such that $\mu_1(B) = \mu (B +1/3)$ for any Borel set $B$ and denote by
$\mu_2$ the probability such that $\mu_2(B) = \mu (B +2/3)$ for any Borel set $B$.

The measure $\mu_3 = \mu_1 + \mu_2.$

Then, $\mu_1$ is absolutely continuous with respect to $\mu_3$. Denote by $R_1$ the Radon-Nikodym derivative.

Moreover, $\mu_2$ is absolutely continuous with respect to $\mu_3$. Denote by $R_2$ the corresponding Radon-Nykodim derivative.

Therefore, for any continuous function $h$ we have
$$ \int h d \mu_1= \int h(z-1/3) d \mu(z)= \int h (z-1/3) R_1(z-1/3) d \mu(z) + \int h (z-2/3) R_1(z-2/3) d \mu(z) $$
and
$$ \int h d \mu_2= \int h(z-2/3) d \mu(z)= \int h (z-1/3) R_2(z-1/3) d \mu + \int h (z-2/3) R_2(z-2/3) d \mu(z). $$

Taking above $h(z) =g(z+1/3)$ the first condition can be rewritten as: for any continuous function $g$:
$$  \int g d \mu= \int g (z) R_1(z-1/3) d \mu(z) + \int g (z-1/3) R_1(z-2/3) d \mu(z) .$$

Taking above $h(z) =g(z+2/3)$ the first condition can be rewritten as: for any continuous function $g$:
$$  \int g d \mu= \int g (z+1/3) R_2(z-1/3) d \mu(z) + \int g (z) R_2(z-2/3) d \mu(z) $$

We will show that
\begin{equation} \label{cio1} \tilde{J}(z) = J(z-2/6) R_1(2 z) + J(z-4/6) R_2(2 z ).\end{equation}

This corresponds also to
\begin{equation} \label{cio4} \tilde{J}(x/2+5/6) = J(x/2 +3/6) R_1(x-1/3) + J(x/2+1/6) R_2(x/2 -1/3 ) \end{equation}
and
\begin{equation} \label{cio5} \tilde{J}(x/2+2/3) = J(x/2+2/6) R_1(x- 2/3) + J(x/2) R_2(x/2 - 2/3 ). \end{equation}

It will follow that the entropy of $\nu$ is
\begin{equation} \label{cio3} h(\nu) = - \int \log  (\,\,J(z-2/6) R_1(2 z) + J(z-4/6) R_2(2 z )\,\,)\, d \nu.\end{equation}

\medskip

\begin{remark}  \label{opt} Note that for any $x$ we have that $R_1(x) + R_2 (x)=1$. The above expression for  $\tilde{J}$ in (\ref{cio1}) says in some sense that $\tilde{J}$ attain values on the convex hull of the values of $J$. It is reasonable to guess that this mechanism is responsible for the increase of entropy under convolution (see a kind of more general and analytic statement in the appendix). Expression (\ref{cio3}) permits an analytic estimation of this increase.
\end{remark}
\medskip

Assuming (\ref{cio1}) we have to show that (\ref{klm34}) is true for any $f$.

\bigskip

a)  Consider first a function $f$ with support on $(0,1/6).$

We have to show that

$$ \int_{2/3}^1 {J(x/2)\,  f(x/2\,+\,2/3\,)}\, \ d\mu(x)  + \int_{1/3}^{2/3} { J(x/2\,+\, 1/2) \, f(\,x/2\,+\,5/6\,)}\,  \, d \mu (x) \,= $$
\begin{equation} \label{a4} \int_{2/3}^1 \,{\tilde{J}\,(x/2\,+\,\,2/3\,\,)  \,f(\,x/2\,+\,\,2/3\, )}\,d \mu(x)  \,
+  \,\int_{1/3}^{2/3}  {\tilde{J}(x/2\,+\,5/6\,) \,f(\,x/2+\,\,5/6\, )}\, d \mu (x) \,\,.\end{equation}

From (\ref{cio1}), (\ref{cio4}) and (\ref{cio5}) we get
$$ \int_{2/3}^1 J(x/2)\,  f(x/2\,+\,2/3\,)\, \ d\mu(x)  + \int_{1/3}^{2/3}  J(x/2\,+\, 1/2) \, f(\,x/2\,+\,5/6\,)\,  \, d \mu (x) \,= $$
$$ [\,\int_{2/3}^1 J(x/2)\,  f(x/2\,+\,2/3\,)\,R_2(x-2/3)  \ d\mu(x) + \int_{1/3}^{2/3} J(x/2+1/6)\,  f(x/2\,+\,5/6\,)\,R_2(x-1/3)  \ d\mu(x)\,] \,+$$
$$ [\int_{2/3}^1 J(x/2+2/6) f(x/2+2/3)R_1(x-2/3)   d\mu(x) + \int_{1/3}^{2/3} J(x/2+3/6)  f(x/2+\,5/6)R_1(x-1/3)   d\mu(x)] =$$
$$[\,\,\int_{2/3}^1 J(x/2) f(x/2+2/3)R_2(x-2/3)   d\mu(x) +\int_{2/3}^1 J(x/2+2/6) f(x/2+2/3)R_1(x-2/3)   d\mu(x)\,\,]+$$
$$[\,\int_{1/3}^{2/3} J(x/2+3/6)  f(x/2+\,5/6)R_1(x-1/3)   d\mu(x)+  \int_{1/3}^{2/3} J(x/2+1/6)\,  f(x/2\,+\,5/6\,)\,R_2(x-1/3)  \ d\mu(x)\,]=$$
$$[\,\,\int_{2/3}^1 f(x/2+2/3)\,\,  (\,J(x/2) R_2(x-2/3)   d\mu(x) + J(x/2+2/6) R_1(x-2/3) \,) \,\, d\mu(x)\,\,]+$$
$$[\,\int_{1/3}^{2/3} f(x/2+\,5/6)\,\,\,(\, J(x/2+3/6)  R_1(x-1/3)   d\mu(x)+  \ J(x/2+1/6)\, \,R_2(x-1/3)  \,) d\mu(x)\,]=$$
$$\int_{2/3}^1 f(x/2+2/3)\,\,  \tilde{J} (x/2+ 2/3) \,\, d\mu(x)\,\,+\,\int_{1/3}^{2/3} f(x/2+\,5/6)\,\, \tilde{J} (x/2+5/6)  d\mu(x)\, ,$$
and this shows (\ref{a4}).

b) Suppose $f$ has support on the interval $[2/6,3/6).$

 We have to show that

$$ \int_{0}^{1/3} {J(x/2)\,  f(x/2\,+\,1/3\,)}\, \ d\mu(x)  + \int_{1/3}^{2/3} { J(y/2+1/2) \, f(\,y/2\,+\,1/6\,)}\,  \, d \mu (y) \,= $$
\begin{equation} \label{a6} \int_{0}^{1/3} \,{\tilde{J}\,(x/2\,+\,\,1/3\,\,)  \,f(\,x/2\,+\,\,1/3\, )}\,d \mu(x)  \,
+  \,\int_{1/3}^{2/3}  {\tilde{J}(y/2 +\,\,1/6) \,f(\,y/2+\,\,1/6\, )}\, d \mu (y) \,\,. \end{equation}

\medskip
The proof is similar to the previous case and it will be left for  the reader.

\medskip

\medskip

c) Suppose the function $f$ has  support on $(4/6,5/6).$  We have to show that

$$ \int_{0}^{1/3} {J(x/2)\,  f(x/2\,+\,2/3\,)}\, \ d\mu(x)  + \int_{2/3}^{1} { J(y/2) \, f(\,y/2\,+\,1/3\,)}\,  \, d \mu (y) \,= $$
\begin{equation} \label{a5}  \int_{0}^{1/3} \,{\tilde{J}\,(x/2\,+\,\,2/3\,\,)  \,f(\,x/2\,+\,\,2/3\, )}\,d \mu(x)  \,
+  \,\int_{2/3}^{1}  {\tilde{J}(y/2\,+\,1/3\,) \,f(\,y/2+\,\,1/3\, )}\, d \mu (y) \,\,. \end{equation}

\medskip
The proof is similar to the previous case and it will be left for  the reader.

For functions $f$ with support on the other possible intervals we proceed in a similar way. This will give the explicit expression of $\tilde{J}$ in terms of $J, R_1, R_2$ on all points.

\subsection{A class of examples where $J=\tilde{J}$} \label{JJ}

\noindent

Suppose we ask: when $J = \tilde{J}$?
What equation should $J$  satisfy in this case? Is there some special form of $J$ such that this happens? We will present examples where this happens.

\medskip

Denote by $\mathcal{S}$ the class of positive H\"older Jacobians $J$ such that for any $x\in [0,1/6)$ we have

\begin{equation} \label{haha} J(x) = J(x- 2/6)= J(x-4/6).
\end{equation}

We point out that under the above conditions the values of $J$ on $[0,1/6)$  determine $J$ uniquely. Indeed, on the intervals $[0,1/6)$ $[2/6,3/6)$ and $[4/6,5/6]$ is clearly determined. On the intervals of the form $[1/6,2/6)$ $[3/6,4/6)$ and $[5/6,1)$ it is also determined because the sum of $J$ on the preimages of any point is equal to $1$.

There exist several continuous (and H\"older) Jacobians satisfying such conditions.

\medskip

\medskip

The equation
$$J(z) = J(z-2/6) R_1(2 z) + J(z-4/6) R_2(2 z )$$
is true for any $f \in \mathcal{S}$ and any $z \in S^1$ because $R_1 + R_2=1$.

As from (\ref{cio1}) we get that  for any $z$
$$\tilde{J}(z) = J(z-2/6) R_1(2 z) + J(z-4/6) R_2(2 z ),$$
it follows that in this case $\tilde{J}= J$ and $\mu=\mu * \rho.$
\medskip

\section{Differentiability of the entropy of convolution} \label{dif}

To each equilibrium probability for a H\"older potential $A$ one can associate a unique positive Holder Jacobian.
Therefore, the set of equilibrium probabilities can be considered as a Banach manifold $\mathcal{N}$ (see \cite{KGLM}). In this way we can consider the bijective map $\log J \to \mu_{\log J}$ over $\mathcal{N}$.

Given a probability $\mu_{\log J}\in \mathcal{N}$ (associated to the potential $\log J$) and a tangent vector $\zeta\in T_{\mu_{\log J} } \mathcal{N} $, one is interested on the derivative $\mu_{ \log + \zeta}$ along $\zeta$, where
$\mu_{ \log J + \zeta}$
is the equilibrium probability for the potential $\log J + \zeta.$

For  a fixed $\varphi$ consider the transformation $G_\varphi$, such that,   $G_\varphi(\log J)=\int\varphi\,d\mu_{\log J}$,
then,
$$D (G_\varphi)_{\log J}(\zeta)
 = \int (I-\mathcal{L}_{\log J})^{-1}(\varphi) \cdot \zeta \,d\mu_{\log J}=$$
 $$
\sum_{i=0}^{+\infty} \int \varphi \cdot \zeta\circ T^i \,d\mu_{\log J}.$$

Given a H\"older potential $A$, following  \cite{KGLM}, denote $\mathcal{N}(A) =\log J$, where $J$ is the Jacobian of the equilibrium probability for $A$. We also denote $\mu_{\mathcal{N}(A)}$ the Gibbs (equilibrium) probability for $A$.

We denote $\mu_i$, $i=1,2$, the probability associated, respectively, to the Jacobians $(\log J)_i.$

We denote $\mu^t =  \mu_{\mathcal{N} (\log J_1 + t z_3)}$, where $z_3$ is a tangent vector to the manifold of Gibbs probabilities at the point $\mu_1$. Note that in this case $\int z_3\,\, d \mu_1=0$.

Denote by $J_1^t$ the Jacobian of $\mu^t$. This means that $\log J_1^t = \mathcal{N} (\log J_1 + t z_3).$

If $\nu_t = \mu^t * \mu_2$ we get that

$$h(\nu_t)=- \int\,[\,\,\int\, \log \,(\,\int J_1^t(r+s-x) d \mu_2(x)\,) \, d  \mu_2(r)\,\, ]\,\,d \mu^t (s)$$

Denote
$$Z^t (s) =\int\, \log \,(\,\int J_1^t(r+s-x) d \mu_2(x)\,) \, d  \mu_2(r)\,\, .$$
Then,
$$ \frac{d}{dt} h ( \nu_t )|_{t=0}= - \int \frac{d}{dt}|_{t=0}\,\, Z^t (s) \, d \nu_t |_{t=0}(s) - \int Z^t (s) \,\, |_{t=0} \,\,\frac{d}{dt}|_{t=0}  d \nu_t(s)  .  $$

Given a continuous function $\phi$ we have from \cite{KGLM} that
$$\int \phi (s) \,\, \,\,\frac{d}{dt}|_{t=0}  d \nu_t(s)  = \int \phi\, z_3\, d \mu_1.  $$

Note that $Z^t - Z^0$ goes uniformly to zero when $t \to 0.$

Therefore, from \cite{KGLM}
$$\int Z^t (s) \,\, |_{t=0} \,\,\frac{d}{dt}|_{t=0}  d \nu_t(s)  =  $$
$$\int Z^0 (s) \,\,  \,\,\frac{d}{dt}|_{t=0}  d \nu_t(s)  =$$
$$\int \,\,[\int\, \log \,(\,\int J_1(r+s-x) d \mu_2(x)\,) \, d  \mu_2(r)]  \,\,\frac{d}{dt}|_{t=0}  d \nu_t(s)= $$
$$ \int \,\,[\int\, \log \,(\,\int J_1(r+s-x) d \mu_2(x)\,) \, d  \mu_2(r)]  \,\,z_3 (s) d \mu_1(s). $$

We denote by $\varphi^t$ and $\lambda^t$, respectively, the main eigenfunction and the main eigenvalue of the
Ruelle operator for the potential  $\log ( J_1) + t z_3.$

Note that when $t=0$ we get that $\varphi^t=1$ and $\lambda^t=1$.

As
$$ \log ( J_1^t) = \log (J_1) +\,t\, z_3\,+ \log \varphi_t - \log(\varphi_t \circ T)- \log \lambda^t,$$
which means
$$ J_1^t = J_1\, e^{t\, z_3 + \log \varphi_t - \log (\varphi_t \circ T)- \log \lambda^t},$$

we get that

$$Z^t (s) =\int\, \log \,(\,\int J_1\, e^{t\, z_3 + \log \varphi_t - \log (\varphi_t \circ T)- \log \lambda^t} (r+s-x) d \mu_2(x)\,) \, d  \mu_2(r)\,\, .$$

and

$$ \frac{d}{dt}|_{t=0}\,\, Z^t (s) =$$
$$\int\,\frac{d}{dt}|_{t=0}\,\, [\,\,\log \,(\,\int J_1\, e^{t\, z_3 + \log \varphi_t - \log (\varphi_t \circ T)- \log \lambda^t} (r+s-x)\,\,\,\, d \mu_2(x)\, ]\, d  \mu_2(r)\,\, .$$

Denote
$$ Y^t (s,r) = \,\,\,(\,\int J_1\,(r+s-x) e^{t\, z_3 + \log \varphi_t - \log (\varphi_t \circ T)- \log \lambda^t} (r+s-x)\,\, d \mu_2(x).$$

 Therefore,
$$  - \int \frac{d}{dt}|_{t=0}\,\, Z^t (s) \, d \nu_t |_{t=0}(s) =  - \int \frac{\frac{d}{dt}|_{t=0}\,\, Y^t (s,r)}{ Y^0 (s,r) } \, d \mu_2(r) \, d \mu_1(s).$$

Now we estimate
$$\frac{d}{dt}|_{t=0}\,\, Y^t (s,r) =  $$
$$\int J_1 \,\,(r+s-x)\,[\,z_3 +  \frac{d}{dt}|_{t=0} \,(\, \log \varphi_t - \log (\varphi_t \circ T)- \log \lambda^t\,)  (r+s-x) \,] \, d \mu_2 (x) .$$

Finally,
$$  - \int \frac{d}{dt}|_{t=0}\,\, Z^t (s) \, d \nu_t |_{t=0}(s) =  $$
$$ -     \int \frac{\int J_1 \,\,(r+s-x)\,[\,z_3 +  \frac{d}{dt}|_{t=0} \,(\, \log \varphi_t - \log (\varphi_t \circ T)- \log \lambda^t\,)   \,] (r+s-x)\, d \mu_2 (x)}{ \int J_1 \,\,(r+s-x)\,   d \mu_2(x)}         \,   d \mu_2(r) \, d \mu_1(s)     $$

In this way we get the following proposition:

\begin{proposition} \label{ert}
Suppose  $\mu_i$, $i=1,2$ are  probabilities associated, respectively, to the Jacobians $\log J_i.$

Denote $\mu^t =  \mu_{\mathcal{N} (\log J_1 + t z_3)}$, $t\in \mathbb{R}$ small,  where $z_3$ is a tangent vector to the manifold of Gibbs probabilities at the point $\mu_1$, and  $\nu_t = \mu^t * \mu_2$.

We also denote by $\varphi^t$ and $\lambda^t$, respectively, the main eigenfunction and the main eigenvalue of the
Ruelle operator for the potential  $\log ( J_1) + t z_3.$

Then,

$$ \frac{d}{dt} h ( \nu_t )|_{t=0}=$$
$$ -     \int \frac{\int J_1 \,\,(r+s-x)\,[\,z_3 +  \frac{d}{dt}|_{t=0} \,(\, \log \varphi_t - \log (\varphi_t \circ T)- \log \lambda^t\,)   \,] (r+s-x)\, d \mu_2 (x)}{ \int J_1 \,\,(r+s-x)\,   d \mu_2(x)}         \,   d \mu_2(r) \, d \mu_1(s)     $$
$$ -\, \int \,\,[\int\, \log \,(\,\int J_1(r+s-x) d \mu_2(x)\,) \, d  \mu_2(r)]  \,\,z_3 (s) d \mu_1(s). $$
\end{proposition}

\section{Appendix} \label{ape}

Now we will prove a result inspired by the reasoning followed  in section \ref{per} (see Remark \ref{opt}).

It a result of interest in itself.

\begin{proposition} \label{por}
Suppose are given the H\"older Jacobians  $J_1$ and $J_2$ and they are such that: $J_2\geq J_1$ when $J_1\geq1/2$, and $J_2\leq J_1$ when $ J_1\leq 1/2$.

Denote $\mu_i$ the Gibbs probability associated to the H\"older potential $\log J_i$, $i=1,2$. Then, $h(\mu_1)\geq h(\mu_2)$.

\end{proposition}

{\bf Proof:}

One way to get a path from $J_1$ to $J_2$ is to take $J^{\,t}= J_1 + t\,(J_2-J_1),$ $t\in [0,1]$.

Note that $J^{\,t} (x_1) + J^t (x_2)=1$, if $T(x_1)=T(x_2)$ (therefore $J^{\,t}$ is a H\"older Jacobian for each value $t$).

We know that if $\int \,\chi\, d \mu_1=0$, then, the entropy $h_t$ of the Gibbs state associated to $\log J_1 +\, t \,\chi$ satisfies

$$ \frac{d\, h_t}{dt}_{t=0}= -\int \,\chi\,\log J_1\, d \mu_1 =  -\int \,\chi\,(\log J_1- \log 1/2)\, d \mu_1$$

(see page 38 in \cite{KGLM}).

In this way if $\chi(x)\geq 0$ when $(\log J_1(x) - \log 1/2)\geq 0$  and $\chi(x)\leq 0$ when $(\log J_1(x) - \log 1/2)\leq 0$  we get that the entropy {\bf decreases} when we go in the direction $\chi$ beginning on $\mu_1$.
This is so because $-\int \,\chi\,\log J_1\,\, d \mu_1<0$.

Take $\epsilon(t)$ such that
$$ \log J_1 + \epsilon(t)= \log (\,J_1 + t\, (J_2-J_1)\,).$$

Note that $ \log J_1 + \epsilon(1)= \log (\,J_2)$.

Then, $\frac{d}{dt} \epsilon (t)|_{t=0}= \frac{J_2}{J_1} -1$.

Moreover,
$$ \int\,( \frac{J_2}{J_1} -1) \, d \mu_1 = \int\, \frac{J_2}{J_1}  \, d \mu_1 -1= $$
$$\int\,\mathcal{L}_{\log J_1}\, (\frac{J_2}{J_1}  \,)\, d \mu_1 -1 = \int\,\mathcal{L}_{\log J_2}\, (1) d \mu_1\,-1=0  $$

\medskip

The proof that $\frac{d}{dt} \epsilon (t)|_{t=0}\leq 0$ is similar to the case $t=0.$

We denote $\mu^t$ the equilibrium state for the normalized potential $\log (J_1) + \epsilon(t).$

Moreover, $\frac{d}{dt} \epsilon (t)|_{t}= \frac{J_2- J_1}{J_1\,-t\, (J_2-J_1)} $.

In this case

$$ \int  \frac{d}{dt} \epsilon (t)|_{t}  \, d \mu^t = \int \frac{J_2- J_1}{J_1\,+t\, (J_2-J_1)} d \mu^t = $$
$$\int \mathcal{L}_{\log ( J_1\,+t\, (J_2-J_1) ) } \,(\,\frac{J_2- J_1}{J_1\,+t\, (J_2-J_1)}\,) d \mu^t = $$
$$ \int \mathcal{L}_0 (J_1-J_2) d \mu^t =0. $$

Then, $\frac{d}{dt} \epsilon (t)|_{t}=  \chi_{\,t} $, $t \in [0,1]$ is  tangent vector on $\mathcal{N}$ at $\log J_{\,t}$.

Moreover,
$$ \frac{d\, h_t}{dt}|_{t}= -\int \,\chi_{\,t}\,\log J_t\, d \mu_t =  -\int \,\chi_{\,t}\,(\log J_t- \log 1/2)\, d \mu_t= $$
 $$ -\int \,\frac{J_2- J_1}{J_1\,+t\, (J_2-J_1)}\,(\log (\,J_1 + t (J_2-J_1)\,) - \log 1/2)\, d \mu_t.  $$

Remember that $J_2\geq J_1$ when $J_1\geq 1/2$, and $J_2\leq J_1$ when $ J_1\leq 1/2$.

When,  $J_2- J_1\geq 0,$ we get that $(\log (\,J_1 + t (J_2-J_1)\,) - \log 1/2)\geq 0.$

On the other hand when  $J_2- J_1\leq 0,$ we get that $(\log (\,J_1 + t (J_2-J_1)\,) - \log 1/2)\leq 0.$

Therefore,
$ \frac{d\, h_t}{dt}|_{t}\leq 0.$

\qed

\end{document}